\documentclass[12pt]{article}

\usepackage{amssymb,amsmath,latexsym}
\usepackage{exscale,latexsym,amsthm,graphics}
\usepackage{epsfig}
\usepackage{psfrag}
\begin{document}

\newtheorem{thm}{Theorem}[section]
\newtheorem{lemma}[thm]{Lemma}
\newtheorem{defn}[thm]{Definition}
\newtheorem{prop}[thm]{Proposition}
\newtheorem{corollary}[thm]{Corollary}
\newtheorem{remark}[thm]{Remark}
\newtheorem{example}[thm]{Example}

\numberwithin{equation}{section}

\def\ee{\varepsilon}
\def\qed{{\hfill $\Box$ \bigskip}}
\def\MM{{\cal M}}
\def\BB{{\cal B}}
\def\LL{{\cal L}}
\def\FF{{\cal F}}
\def\EE{{\cal E}}
\def\QQ{{\cal Q}}
\def\AA{{\cal A}}
\def\CC{{\cal C}}
\def\q{{\cal q}}
\def\R{{\bf R}}
\def\N{{\mathbb N}}
\def\E{{\bf E}}
\def\F{{\bf F}}
\def\H{{\bf H}}
\def\P{{\bf P}}
\def\Q{{\bf Q}}
\def\S{{\bf S}}
\def\J{{\bf J}}
\def\K{{\bf K}}
\def\F{{\bf F}}
\def\A{{\bf A}}
\def\loc{{\bf loc}}
\def\eps{\varepsilon}
\def\semi{{\bf semi}}
\def\wh{\widehat}
\def\pf{\noindent{\bf Proof.} }
\def\dim{{\rm dim}}

\title{\Large \bf Intrinsic Ultracontractivity for Non-symmetric
L\'evy Processes}
\author{Panki Kim\\
Department of Mathematics\\
Seoul National University\\
Seoul 151-742, Republic of Korea\\
Email: pkim@snu.ac.kr \smallskip \\
and
\smallskip \\
Renming Song\thanks{The research of this author is supported in part
by a joint
US-Croatia grant INT 0302167.}\\
Department of Mathematics\\
University of Illinois \\
Urbana, IL 61801, USA\\
Email: rsong@math.uiuc.edu }

\maketitle

\begin{abstract} Recently in \cite{KS4, KS5}, we extended the concept of
intrinsic ultracontractivity to non-symmetric semigroups and proved
that for a large class of non-symmetric diffusions $Z$ with
measure-valued drift and potential, the semigroup of $Z^D$ (the
process obtained by killing $Z$ upon exiting $D$) in a bounded
domain is intrinsic ultracontractive under very mild assumptions.

In this paper, we study the intrinsic ultracontractivity for
non-symmetric discontinuous L\'{e}vy processes. We prove that, for a
large class of non-symmetric discontinuous L\'{e}vy processes $X$
such that the Lebesgue measure is absolutely continuous with respect
to the L\'{e}vy measure of $X$, the semigroup of $X^D$ in any
bounded open set $D$ is intrinsic ultracontractive. In particular,
for the non-symmetric stable process $X$ discussed in \cite{V}, the
semigroup of $X^D$ is intrinsic ultracontractive for any bounded set
$D$. Using the intrinsic ultracontractivity, we show that the
parabolic boundary Harnack principle is true for those processes.
Moreover, we get that the supremum of the expected conditional
lifetimes in a bounded open set is finite. We also have results of
the same nature when the L\'{e}vy measure is compactly supported.
\end{abstract}

\vspace{.2truein}

\bigskip

\bigskip

\begin{doublespace}

\section{Introduction}
Suppose that $H$ is a semi-bounded self-adjoint operator on $L^2(D)$
with $D$ being an open set in $\R^d$ and that $\{e^{Ht}\}$ is an
irreducible positivity-preserving semigroup with integral kernel
$a(t, x, y)$. We assume that the top of the spectrum $\lambda_1$ of
$H$ is an eigenvalue. In this case, $\lambda_1$ has multiplicity one
and the corresponding eigenfunction $\phi_1$, normalized by
$\|\phi_1\|_{L^2(D)}=1$, is positive almost everywhere on $D$.
$\{e^{Ht}\}$ is said to be intrinsic ultracontractive if for every $
t>0$, there exists $c_t \in (0, \infty)$ such that $a(t, x, y) \le
c_t \phi_1 (x) \phi_1 (y)$.

The notion of the intrinsic ultracontractivity above was introduced
in \cite{DS}. It is a very important concept in both analysis and
probability, and has been studied extensively. When $H$ is the
Dirichlet Laplacian in a domain $D$ (equivalently, the corresponding
process is a killed Brownian motion), the semigroup  $\{e^{Ht}\}$ is
intrinsic ultracontractive for a large class of non-smooth domains
(see, for instance \cite{Ba, BB}). For symmetric $\alpha$-stable
processes with $\alpha\in (0, 2)$, the intrinsic ultracontractivity
has been discussed in \cite{CS1, CS2, Ku}. After obtaining the main
results of this paper, we found out from \cite{GR} that the
intrinsic ultracontractivity for some large classes of symmetric
L\'evy processes was studied in \cite{G}.

Very recently in \cite{KS4}, we extended the concept of intrinsic
ultracontractivity to non-symmetric semigroups and, by using an
analytic method, we proved there that the semigroup of a killed
diffusion process in a bounded Lipschitz domain is intrinsic
ultracontractive if the coefficients of the generator of the
diffusion process are smooth. In \cite{KS5}, by using a
probabilistic  method we proved that for a non-symmetric diffusion
with measure-valued drift and potential belonging to appropriate
Kato classes, the semigroup of the killed process in a bounded
domain is intrinsic ultracontractive when the bounded domain is one
of  the following types: twisted H\"{o}lder domains of order $\alpha
\in (1/3, 1]$, uniformly H\"{o}lder domains of order $\alpha \in (0,
2)$ and domains which can be locally represented as the region above
the graph of a function (see \cite{KS5} for details).

In this paper, we continue our discussion of intrinsic
ultracontractivity for non-symmetric semigroups. We study  the
intrinsic ultracontractivity for non-symmetric discontinuous
L\'{e}vy processes under one of the following two non-overlapping
assumptions on the L\'evy measure: the first case is that the
Lebesgue measure is absolutely continuous with respect to the L\'evy
measure and the Radon-Nikodym derivative is locally integrable away
from $0$ and the second case is that the L\'evy measure is compactly
supported. In the first case, we show that for any bounded open set,
the semigroup of the killed process is intrinsic ultracontractive if
the transition density of the killed process is strictly positive,
bounded and continuous. In particular, the semigroup of the killed
strictly $\alpha$-stable process in any bounded open set is
intrinsic ultracontractive. In the second case we put some mild
assumptions on both the open set and the L\'evy measure:  We assume
that the open set is bounded $\kappa$-fat (a disconnected analogue
of John domain, for the definition see Definition \ref{fat}) and
that the Radon-Nikodym derivative of the absolutely continuous part
of L\'evy measure is bounded below by a positive constant near the
origin. We show that in this case, the intrinsic ultracontractivity
is true if the transition density of the killed process is strictly
positive, bounded and continuous. We do not assume that our
non-symmetric L\'{e}vy process is a purely discontinuous process. It
may contain diffusion and drift parts.

The content of this paper is organized as follows. In Section 2, we
recall some preliminary facts about non-symmetric L\'evy processes.
Section 3 contains the proof of the intrinsic ultracontractivity. We
also show in Section 3 that the intrinsic ultracontractivity implies
the parabolic boundary Harnack principle and that the supremum of
the expected conditional lifetimes is finite. In the last section we
collect some concrete examples of non-symmetric L\'evy processes
satisfying the assumptions of this paper.

In this paper
we use the convention $f(\partial)=0$.
In this paper we will also use the following convention:
the values of the constants $c_1, c_2, \cdots$
might change from one appearance to another. The labeling
of the constants $c_1, c_2, \cdots$ starts anew in the
statement of each result.

In this paper, we use ``$:=$" to denote a
definition, which is  read as ``is defined to be".

\section{Non-symmetric L\'{e}vy Processes}

Let $X=(X_t, \P_x)$ be a L\'{e}vy process in $\R^d$ with the generating
 triplet $(A, \nu, \gamma)$.
i.e., for every
$z\in \R^d$,
$$
\E_0\left[e^{iz\cdot X_1}\right] =\exp \left( -\frac{1}{2}z\cdot Az
+ i\gamma \cdot z +\int_{\R^{d}} (e^{iz \cdot x}-1-iz \cdot x 1_{\{
|x|\le 1\}}(x)) \nu(dx)\right)
$$
where $A$ is a symmetric nonnegative definite $d \times d$ matrix,
$\gamma \in \R^d$, and
$\nu$ is a measure on $\R^d$ satisfying
\begin{equation}\label{e:Lm}
\nu(\{ 0 \})=0 \quad \mbox{ and } \quad \int_{\R^d} (|x|^2 \wedge 1)
\nu(dx) < \infty.
\end{equation}
$\gamma$ is called the drift of $X$ and
$\nu$ is called the  L\'{e}vy measure of $X$.

$-X$ is also a  L\'{e}vy process and it is the dual of $X$. For this
reason we sometimes use $\wh X$ to denote this process. From the
above definition, it is clear that $\wh X$ is a L\'{e}vy process in
$\R^d$ with the generating triplet $(A, \nu(-dx), -\gamma)$.

Let
$$
P_tf(x):= \E_x[f(X_t)]\quad \mbox{ and }\quad
\wh P_tf(x):= \E_x[f(\wh X_t)].
$$
Then for any non-negative Borel functions $f$ and $g$,
$$
\int_{\R^d} P_tf(x) g(x)dx =\int_{\R^d}f(x)  \wh P_tg(x)dx.
$$

Throughout this paper, we assume the following.
\begin{description}
\item{(A1)}  The  L\'{e}vy measure  $\nu$ satisfies either (a) or (b)
below:
\begin{description}
\item{(a)} The Lebesgue measure in $\R^d$ is absolutely continuous
with respect to $\nu$. i.e., there exists
non-negative Borel function $L(x)$ such that for any Borel set $B$,
\begin{equation}\label{e:A1a}
|B|= \int_{B} L(x) \nu(dx).
\end{equation}
Moreover, we assume that $L$ is locally integrable on
$\R^d\setminus \{0\}$ with respect to
the Lebesgue measure.
\item{(b)}
Let $M(x)$ be the Radon-Nikodym derivative of the absolutely
continuous part of $\nu$. We assume that
there exists $R_0 >0$ such that
\begin{equation}\label{e:A1b}
\inf_{x \in B(0,R_0)}M(x) >0.
\end{equation}
\end{description}
\end{description}
\medskip

In \cite{KS5}, we have already discussed the case when $\nu=0$.
The second assumption in (a) is the same as assuming
that $L$ is locally $L^2$-integrable on $\R^d\setminus \{0\}$ with respect to
$\nu$.  (b)  covers the case where the L\'{e}vy measures have
compact supports.

Let $C_0(\R^d)$ be the class of bounded
continuous functions $f$ on $\R^d$
with $\lim_{|x|\to\infty}f(x)=0$.
We say a Markov process $Y$ in $\R^d$ has the Feller property if for every
$g \in C_0(\R^d)$, $\E_x[g(Y_t)]$ is in $C_0(\R^d)$.
Any L\'evy process in  $\R^d$
has the Feller property (for example, see \cite{Be, Sa}).

For any open set $U$, we use $\tau_U$ to denote the first exit
time of $U$ for $X$.
i.e., $\tau_U:=\inf\{t>0: \, X_t\notin U\}$.
Given  an open set $U\subset \R^d$, we define
$X^U_t(\omega)=X_t(\omega)$ if $t< \tau_U(\omega)$
and $X^U_t(\omega)=\partial$ if $t\geq  \tau_U(\omega)$,
where $\partial$ is a cemetery state. The process
$X^U$ is called a killed process in $U$.
We  use $\wh \tau_U$ to denote the first exit
time of $U$ for $\wh X$.
i.e., $\wh \tau_U:=\inf\{t>0: \,\wh X_t\notin U\}$.
We similarly define
$\wh X^U$.

For any $t>0$, define
$$
P^D_tf(x):=\E_x [f(X^D_t)]\quad
\mbox{and}\quad
\wh{P}^D_tf(x):=\E_x [f(\wh X^D_t)].
$$
The next equality is known as Hunt's switching identity
(for example, see Theorem
II.5 in \cite{Be}).
$$
\int_{D}f(x)P^D_tg(x)dx=\int_Dg(x)\wh{P}^D_tf(x) dx.
$$
For the remainder of this section, $D$ is a fixed bounded open set
in $\R^d$. The next assumption is needed to define intrinsic
ultracontractivity for non-symmetric semigroups (see \cite{KS4}).

\begin{description}
\item{(A2)}
The transition density function $p^D(t, x, y)$ for $X^D_t$ exists.
Moreover each $t>0$, $p^D(t, \cdot, \cdot)$ is continuous in $ D
\times D$.
\end{description}

We further assume that $p^D(t, \cdot, \cdot)$ is bounded.
\begin{description}
\item{(A3)}
$\{P^D_t\}$ is ultracontractive. i.e., for  $t>0$, there exists
positive constant $c_t$ such that
$$p^D(t, x, y) \, \le \, c_t \,<\, \infty, \quad (x,y) \in D \times D.
$$
\end{description}

\begin{remark}\label{rk1} {\rm
We do not  know any necessary and sufficient conditions for
(A2)-(A3) in terms of the L\'evy measure. In fact, no necessary and
sufficient condition in terms of the L\'evy measure for the
existence of transition density for L\'evy process is known (see
\cite{Sa} for some sufficient conditions). }
\end{remark}

\medskip

In the remainder of this section, we discuss some elementary
consequences of (A2)-(A3).

From Hunt's switching identity and the continuity of $p^D(t,x,y)$,
we see that $\wh p^D(t, x, y):= p^D(t, y, x)$  is the transition
density for $\wh X$.
 So for every $t>0$ and
Borel set $A \subset D$,
\begin{equation}\label{e:kd1}
\P_x(X^D_t \in A) \, =\,\int_{A} p^D(t,x,y) dy  \quad \mbox{ and } \quad
\P_x(\wh X^D_t \in A) \, =\,\int_{A} p^D(t,y,x) dy.
\end{equation}

If $U \subset D$, then for every $t>0$, $x \in U$ and nonnegative
Borel function $f$,
$$
P^U_tf(x) \le \int_U p^D(t,x,y) f(y)dy \le c_t \int_U  f(y)dy.
$$
Thus $\P_x(X_t^U \in dy)$ is absolutely continuous with respect to
the Lebesgue measure and for every $t>0$, $x \in U$ the density
$p^U(t,x, \cdot)$ exists. Similarly, if $U \subset D$, $\P_x(\wh
X_t^U \in dy)$ is absolutely continuous with respect to the Lebesgue
measure and for every $t>0$, $x \in D$ the density $\wh p^U(t,x,
\cdot)$ exists. Moreover, from Hunt's switching identity, we see
that for every $t>0$,
$$
p^U(t, x, y) =\wh p^U( t, y,x ), \quad \mbox{a.e.}~ (x,y) \in D \times D.
$$
Note that in general we do not know whether $p^U(t,x,y)$ and
$\wh p^U( t, y,x )$ are continuous and strictly positive.

From Lemma 48.3 in \cite{Sa},  it is easy to see that for any bounded open
subset $U$,
there exists $t_1 >0$ such that $\sup_{x \in \R^d} \P_x(X_{t_1} \in U) <1$.
Thus
$$
\theta :=\sup_{x \in \R^d} \P_x(\tau_U > t_1 )\, \le \, \sup_{x \in
\R^d} \P_x(X_{t_1} \in U)\, <\,1.$$ By the Markov property and an
induction argument,
$$ \sup_{x \in \R^d} \P_x(\tau_U > n t_1 ) \,\le\, \theta^n.$$
Thus
\begin{equation}\label{e:ebd}
\sup_{x \in U} \E_x[\tau_U] \,\le\, \frac{t_1}{1-\theta} \,<\, \infty
\end{equation}
(see \cite{C2} for the details).

For any bounded  open subset $U \subset D$, we will use
$G_U(x,y)$ to denote the Green function of $X$ in $U$. i.e.,
$$
G_U(x,y) := \int_0^{\infty} p^U(t,x,y)dt, \quad (x,y) \in U\times U.
$$
By (\ref{e:ebd}),
 $$
\E_x[\tau_U] = \int_U G_U(x,y) dy \,<\, \infty, \quad x \in U
$$
and
$G_U(x,\cdot)$ is well-defined a.e. $U$.

Also  (\ref{e:ebd}) implies that
for every open set $U\subset D$ and $A \subset U^c$ with dist$(A,U)>0$,
we have
\begin{equation}\label{e:P_f1}
\P_x\left(X_{\tau_U} \in A \right)
=
\int_{U}
G_U(x,y) \nu(y-A) dy.
\end{equation}
(for example, see \cite{IW}).

Similarly the Green function $\wh G_U(x,y)$ of $X$ in $U$ is defined as
$$
\wh G_U(x,y) := \int_0^{\infty} \wh p^U(t,y,x)dt, \quad (x,y) \in
U\times U,
$$
which is well-defined a.e. $U$.
For every $A \subset U^c$ with dist$(A,U)>0$, we have
 \begin{equation}\label{e:P_f2}
\P_x\left(\wh X_{\tau_U} \in A \right)
=\int_{U}
\wh G_U(x,y) \nu(A-y) dy.
\end{equation}
Clearly,
$$
G_U( x, y) =\wh G_U( y,x ), \quad \mbox{a.e.}~ (x,y) \in U \times U.
$$

\section{Intrinsic Ultracontractivity for Non-symmetric L\'{e}vy Processes}

In this section, we first recall the definition of the intrinsic
ultracontractivity for non-symmetric semigroups from \cite{KS4} and
then prove that the intrinsic ultracontractivity is true if the
killed non-symmetric L\'evy process $X^D$ satisfies (A1)-(A3) in the
previous section and (A4)-(A5) below. We will use some ideas from
\cite{Ku}.

Many results in this section are stated for both $X^D$ and its dual
$\wh X^D$. Since the proofs for the two processes are similar, we
only present the proofs for $X$.

The following definition is taken from \cite{SW}.

\medskip

\begin{defn}\label{fat} Let $\kappa \in (0,1/2]$. We say that an open set
$D$ in $\R^d$ is $\kappa$-fat if there exists $R>0$ such that for
each $Q \in \partial D$ and $r \in (0, R]$, $D \cap B(Q,r)$ contains
a ball $B(A_r(Q),\kappa r)$ for some $A_r(Q) \in D$. The pair $(R,
\kappa)$ is called the characteristics of the $\kappa$-fat open set
$D$.
\end{defn}

\medskip

Note that every Lipschitz domain and every non-tangentially
accessible domain (see \cite{JK} for the definition of
non-tangentially accessible domains) are $\kappa$-fat. Moreover,
every John domain is $\kappa$-fat (see Lemma 6.3 in \cite{MV}). The
boundary of a $\kappa$-fat open set can be highly nonrectifiable
and, in general, no regularity of its boundary can be inferred.
Bounded $\kappa$-fat open sets may be disconnected.

Depending on whether (A1)(a) or (A1)(b) is valid, our assumptions on
the open set $D$ are different. In both cases, we will need to
define some subsets $B_0$, $C_1$ and $B_2$ of $D$. The following
assumptions on $D$ will always be in force in the reminder of this
section.
\begin{description}
\item{(A4)(a)}
If $\nu$ satisfies (A1)(a),  we assume that $D$ is an
arbitrary bounded open set.
Choose a point $x_0$ in D and $r_0 \in (0, \infty)$ such that
$B(x_0, 2r_0) \subset \overline{B(x_0, 2r_0)} \subset D$.
We put $B_0:= B(x_0, r_0/2)$, $C_1:= \overline{B(x_0, r_0)}$ and
 $B_2:= B(x_0, 2r_0)$.

\item{(A4)(b)}
If $\nu$ satisfies (A1)(b), then we assume that $D$ is a bounded
$\kappa$-fat open with  the characteristics
$(R, \kappa)$.
Without loss of generality, we assume $R \le \frac12 R_0$ where
$R_0$ is the constant in (A1)(b). Let $\rho(x)$ be the
distance of a point $x$ to the boundary of $D$, i.e.,
$\rho(x)=\text{dist} (x,\partial D)$.
Define
\begin{eqnarray*}
B_0&:=&\{x \in D: \rho(x)  > R \kappa /2 \},\\
C_1&:=&\{x \in D: \rho(x)  \ge R \kappa /4 \},\\
B_2&:=&\{x \in D: \rho(x)  > R \kappa /8 \}.
\end{eqnarray*}
\end{description}
The distinction between (A4)(a) and (A4)(b) will be made only in the
proof of Lemma \ref{l:IU1} below.

Define
$$
\eta_U\,:=\,\inf\{t\ge 0: \, X_t\notin U\}
\quad
\mbox{and}
\quad
\wh\eta_U\,:=\,\inf\{t\ge 0: \,\wh X_t\notin U\}.
$$
Note that $\eta_U \le \tau_U$ and $\wh \eta_U \le \wh \tau_U$. Moreover,
$\eta_U(\omega) = \tau_U(\omega)$ and $\wh \eta_U(\omega) =
\wh \tau_U(\omega)$
if $X_0(\omega) \in U$ and $\wh X_0(\omega) \in U$ respectively.

\medskip

\begin{lemma}\label{l:IU1}
If (A1)-(A4) are true, then there exists a constant $c>0 $ such that
for every $x \in \R^d\setminus C_1$,
$$
\P_x \left(X_{\eta_{D\setminus C_1}} \in C_1\right)
\,\ge \, c\,\E_x \left[  \eta_{D\setminus C_1} \right]
\quad
\mbox{and}
\quad
\P_x \left(\wh X_{\wh \eta_{D\setminus C_1}} \in C_1\right)
\,\ge \,c \,\E_x \left[ \wh \eta_{D\setminus C_1} \right].
$$
\end{lemma}

\pf
If $x \in \R^d\setminus D$, $\P_x (\eta_{D\setminus C_1}=0)=1$. Thus $\E_x
[  \eta_{D\setminus C_1} ]=0$ and assertions of the lemma are trivial in
this case.
Now we assume $x \in D \setminus C_1$.
\begin{description}
\item{(1)}
First we deal with the case that $\nu$ satisfies (A1)(a). If $w \in
B_0$ and $y \in D \setminus C_1$, then $|w-y|\ge |y-x_0|-|w-x_0| >
r_0/2$ and $|w-y|<2$diam$(D)$. So the set
$$
A:=\bigcup_{y \in D\setminus C_1} (y-B_0)$$
is a relatively compact subset of $\R^d \setminus \{0\}$.
By (A1), for every $y \in D \setminus C_1$ we have
$$
|B_0|\, =\,|y-B_0| \,\le\,\int_{A} 1_{y-B_0}(z) L(z)\nu(dz) \,\le\,
(\nu(y-B_0))^{1/2} \|1_{A}L \|_{L^2(\nu)}.
$$
We know from our assumption (A1) that
$$
\|1_{A}L \|_{L^2(\nu)}^2 \, =\, \int_{A}L^2(z)\nu(dz) \, =\,
\int_{A}L(z)dz \,<\, \infty.
$$
Therefore from (\ref{e:P_f1}), we have
\begin{eqnarray*}
&&\P_x \left(X_{\eta_{D\setminus C_1}} \in C_1\right)
\ge
\P_x \left(X_{\tau_{D\setminus C_1}} \in B_0 \right)\\
&&=\int_{D\setminus C_1 }
G_{D\setminus C_1}(x,y) \int_{y-B_0} \nu(dz) dy
\ge
\int_{D\setminus C_1 }
G_{D\setminus C_1}(x,y)      \frac{|B_0|^2}{ \|1_{A}L
\|_{L^2(\nu)}^2}  dy\\
&&=   \frac{|B_0|^2}{ \|1_{A}L
\|_{L^2(\nu)}^2} \,
\E_x \left[  \tau_{D\setminus C_1} \right]\,=\, \frac{|B_0|^2}{ \|1_{A}L
\|_{L^2(\nu)}^2} \,  \E_x \left[  \eta_{D\setminus C_1} \right].
\end{eqnarray*}
\item{(2)}
Now we deal with the case that $\nu$ satisfies (A1)(b).
For each $y \in  D \setminus C_1$, choose a point $Q_y \in \partial D$
such that $\rho(y)=|y-Q|< \kappa R /4$. Since $D$ is $\kappa$-fat,
there exists  a point $A_y \in D$
such that $B(A_y,  \kappa R) \subset D \cap B(Q_y, R)$.
It is easy to see that
\begin{equation}\label{e:ka}
B(A_y, \frac12 \kappa R) \subset B_0 \cap  B(Q_y, R)
\subset B_0 \cap  B(y, R_0).
\end{equation}
In fact, if $|w-A_y| <  \frac12 \kappa R$, then
$\rho(w) \ge \rho(A_y) -|w-A_y| >\kappa R-   \frac12 \kappa R
= \frac12 \kappa R$.
If $|w-Q_y| < R$, then $|y-w|\le |y-Q_y|+|w-Q_y| < R +
\kappa R /4 < 2R \le R_0$.
Thus by (\ref{e:ka}) and (A1)(2), we have for every $y \in  D \setminus C_1$
\begin{eqnarray*}
&&\nu(y-B_0) \ge  \int_{B_0} M(y-z)    dz
\ge \int_{ B(A_y, \frac12 \kappa R)  } M(y-z)    dz\\
&&\ge \left(\inf_{w\in B(0, R_0)}M(w) \right)\, | B(0, \frac12
\kappa R)| \,=:\,c_1(R, R_0, \kappa, d) \,>\, 0.
\end{eqnarray*}
Now by (\ref{e:P_f1}), we get
\begin{eqnarray*}
&&\P_x \left(X_{\eta_{D\setminus C_1}} \in C_1\right)
\ge
\P_x \left(X_{\tau_{D\setminus C_1}} \in B_0 \right)\\
&&=\int_{D\setminus C_1 }
G_{D\setminus C_1}(x,y) \nu(y-B_0) dy
\ge  c_1  \,
\E_x \left[  \tau_{D\setminus C_1} \right]\,=\,c_1\,
\E_x \left[  \eta_{D\setminus C_1} \right].
\end{eqnarray*}
\end{description}
\qed

\medskip

Let $\theta$ be the usual shift operator for Markov processes, and we
define stopping times $S_n$ and $T_n$ recursively by
$$
S_1 \,:=\, 0, \quad
T_n\,:=\, S_n +   \eta_{D\setminus C_1} \circ \theta_{S_n}
\quad \mbox{and}\quad
S_{n+1}\,:=\, T_{n} +   \eta_{B_2} \circ \theta_{T_{n}}, \quad n \ge 1.
$$
Similarly we define $\wh T_n$ and $\wh S_n$ for $\wh X$.

\medskip

\begin{lemma}\label{l:IU3} If (A1)-(A4) are true, then there exists a constant $c>0$ such that for every $x \in D$
$$
\P_x (X_{T_n} \in C_1)
\,\ge \, c\,\E_x [ T_n-S_n  ]
\quad
\mbox{ and }
\quad
\P_x (\wh X_{\wh T_n} \in C_1)
\,\ge \, c\,\E_x [\wh T_n-\wh  S_n  ].
$$
\end{lemma}

\pf Since
$T_n= S_n +   \eta_{D\setminus C_1} \circ \theta_{S_n}$, by the
strong Markov property,
$$
\P_x \left(X_{T_n} \in C_1\right)
\,=\,
\P_x \left(X_{S_n +   \eta_{D\setminus C_1} \circ \theta_{S_n}} \in C_1\right)
\,=\, \E_x \left[\P_{X_{S_n}}   \left(\eta_{D\setminus C_1}\in C_1
\right) \right].
$$
Applying Lemma \ref{l:IU1} to the equation above, we get
$$
\P_x \left(X_{T_n} \in C_1\right) \,\ge\, c\, \E_x \left[
\E_{X_{S_n}}\left[  \eta_{D\setminus C_1} \right] \right]\,=\,
c\, \E_x [\eta_{D\setminus C_1} \circ \theta_{S_n} ]
\,=\, c\, \E_x \left[T_n-S_n \right].
$$
\qed

\medskip

\begin{lemma}\label{l:IU6}
For every $x \in D$,
$$
\P_x\left(\lim_{n \to \infty}S_n=\lim_{n \to \infty} T_n
= \tau_D\right)\,=\,\P_x\left(
\lim_{n \to \infty} \wh S_n=\lim_{n \to \infty}\wh T_n
= \wh \tau_D\right)\,=\,1.
$$
\end{lemma}

\pf Recall that  $ \P_x( \tau_D < \infty) =1 $, Clearly $S_n \le
\tau_D$, Let $S:=\lim_{n \to \infty} S_n  \le \tau_D$. we define a
subprocess $Z$ of $X^D$ by letting $Z_t(\omega)=X_t(\omega)$ if $t<
S(\omega)$ and $Z_t(\omega)=\partial$ if $t\geq  S(\omega)$. By
Corollary III.3.16 in \cite{BG}, $Z$ is a Hunt process. Thus by the
quasi-left continuity,
\begin{eqnarray*}
&&\P_x( S< \tau_D)=\P_x( S< \tau_D,  \,\tau_D < \infty )
\,=\,\P_x( Z_{S-} \in D,  \,\tau_D < \infty)\\
&&=\,
 \P_x(\lim_{n \to \infty} Z_{S_n} \in D, S_n < S_{n+1},  \, \tau_D < \infty)
\,\le\, \P_x( Z_{S} \in D,  \,S < \infty) =0.
\end{eqnarray*}\qed

\medskip

By the separation property for Feller processes, there exists $t_0$ such that
\begin{equation}\label{e:FS}
\inf_{y \in C_1} \P_y( \tau_{B_2} > t) \, \ge  \,\frac12 \quad
\mbox{ and } \quad \inf_{y \in C_1} \P_y(\wh \tau_{B_2} > t) \, \ge
\, \frac12
\end{equation}
for any $t\le t_0$ (see Exercise 2 on page 73 of \cite{CW}).

\medskip

\begin{lemma}\label{l:main} If (A1)-(A4) are true, then there exists $c>0$ such that for any $t\le t_0$,
$$
\P_x(X_{t} \in B_2,  \,\tau_D > t)\,\ge \, c \,
 \int_{D\setminus B_2} G_D(x,y)dy
, \quad x \in  D
$$
and
$$
\P_x(\wh X_{t} \in B_2,  \,\wh \tau_D > t)\,\ge \, c \,
 \int_{D\setminus B_2} \wh G_D(x,y)dy
, \quad x \in  D .
$$
\end{lemma}

\pf Note that by Lemma \ref{l:IU6},
\begin{eqnarray*}
\P_x(X_{t} \in B_2,  \, \tau_D > t) &= &\P_x(\cup_{n=1}^{\infty} \{
X_{t} \in B_2,  \,
S_n \le t< S_{n+1}\} )\\
&\ge &\P_x(
\cup_{n=1}^{\infty} \{ X_{t} \in B_2,  \,  T_n \le t< S_{n+1}\})\\
&=& \sum_{n=1}^{\infty} \P_x(X_{t} \in B_2, \,  T_n \le t < S_{n+1})
\,= \,\sum_{n=1}^{\infty} \P_x(T_n \le t < S_{n+1}).
\end{eqnarray*}
By the strong Markov property and (\ref{e:FS}),
\begin{eqnarray*}
&& \P_x(T_n \le t < S_{n+1}) \,=\,
  \P_x(T_n \le t <  T_{n} +   \eta_{B_2} \circ \theta_{T_{n}} )\\
&&= \E_x \left[ \P_{X_{T_n}} (  t  < \eta_{B_2}   )     \right]
\,\ge\,\E_x \left[ \P_{X_{T_n}} (  t  < \tau_{B_2}   ) : X_{T_n} \in
C_1   \right] \,>\, \frac12 \P_x(  X_{T_n} \in C_1 ),
\end{eqnarray*}
which is larger than $ c_1 \E_x[ T_n-S_n]$ for some constant $c_1>0$
by Lemma \ref{l:IU3}.
Therefore by Lemma \ref{l:IU6} and Fubini's theorem, for $x \in D$,
\begin{eqnarray*}
&&\P_x(X_{t} \in B_2, \tau_D > t) \ge  c_1\sum_{n=1}^{\infty} \E_x[
T_n-S_n] =  c_1\sum_{n=1}^{\infty}  \E_x\left [ \int_{S_n}^{T_n}
1_{\R^d}(X_t)
dt\right]\\
&&\ge
c_1 \E_x\left [\sum_{n=1}^{\infty}  \int_{S_n}^{T_n}
1_{D \setminus B_2}(X_t) dt\right]
=\E_x\left [\sum_{n=1}^{\infty} \int_{S_n}^{T_n} 1_{D \setminus B_2}(X_t) dt
+ \sum_{n=1}^{\infty} \int_{T_n}^{S_{n+1}} 1_{D \setminus B_2}(X_t)
dt \right]\\
&&=  c_1\E_x \int_0^{\tau_D} 1_{D \setminus B_2}(X_t) dt\,=\,
c_1 \int_{D \setminus B_2} G_D (x,y)dy.
\end{eqnarray*}\qed

\medskip

The above lemma will also be used in the next section to prove the
strict positivity of the density of killed processes for some
particular non-symmetric L\'evy processes.

 The next proposition is elementary and should
be well-known. But we could not find any reference for this. We
include a proof here for completeness.

\medskip

\begin{prop}\label{SC} For any open set $D$ with finite Lebesgue
measure, $\{P^D_t\}$ and $\{\wh{P}^D_t\}$ are both strongly
continuous contraction semigroups in $L^2(D, dx)$.
\end{prop}

\pf The contraction property follows easily from the duality and
H\"{o}lder's inequality. So we only prove the strong continuity.

Recall that for any open subset $U$ of $\R^d$ and any $x\in U$, we
have
\begin{equation}\label{e:tzero}
\lim_{t \to 0} \P_x(\tau_U \le t)\,=\, \P_x(\tau_U =0) \,=\,0.
\end{equation}
We first consider $f$ in $C_c(D)$, the class of continuous functions
on $D$ with compact supports. Fix $x \in D$. Given $\eps>0$, choose
$\delta >0$ such that $|f(y)-f(x)| < \eps/2$ for $y \in B(x,\delta)
\subset D$. Then for $x \in D$,
\begin{eqnarray*}
&&|P^D_tf(x)-f(x)|\\
&\le & \int_D p^D(t,x,y)|f(y)-f(x)|dy  \,+ \, |f(x)|\P_x(\tau_D \le t)\\
&\le & \left(\int_{D\cap \{ |x-y| < \delta \}} + \int_{D\cap \{
|x-y| \ge \delta \}} \right) p^D(t,x,y)|f(y)-f(x)|dy + \|f\|_\infty
\P_x(\tau_D \le t)\\
&\le &\frac{\eps}2 \,+ \,2 \|f\|_\infty \P_x(|X^D_t-x| \ge \delta) +
\|f\|_\infty
\P_x(\tau_D \le t)\\
&\le &\frac{\eps}2 \,+ \,2 \|f\|_\infty   \P_x(\tau_{B(x, \delta)}
\le t)   + \|f\|_\infty \P_x(\tau_D \le t).
\end{eqnarray*}
Applying (\ref{e:tzero}) to both $\P_x(\tau_{B(x, \delta)} \le t)$
and $\P_x(\tau_D \le t)$, we get that $P^D_tf$ converges pointwise
to $f$. Since $\|  P_t^Df  \|_{\infty} \le\| f  \|_{\infty}$ and $D$
has finite Lebesgue measure, by the bounded convergence theorem,
$P^D_tf$ also converges to $f$ in $L^2(D)$.

Now we assume $f \in L^2(D)$. Given $\eps > 0$, choose $g \in
C_c(D)$ with $\|f-g\|_{L^2(D)} < \eps/4$. By  the contraction
property of $P^D_t$,
\begin{eqnarray*}
\| P^D_tf -f \|_{L^2(D)} &\le& \| P^D_t(f-g) \|_{L^2(D)}
+ \| P^D_tg -g \|_{L^2(D)} + \| f-g \|_{L^2(D)} \\
&\le& 2  \| f-g \|_{L^2(D)} + \| P^D_tg -g \|_{L^2(D)} \,\le\,
\frac{\eps}{2} \,+\, \| P^D_tg -g \|_{L^2(D)} .
\end{eqnarray*}
Thus $P^D_tf$ converges to $f$ in $L^2(D)$. \qed

\medskip

Our last assumption below will be used to define the intrinsic
ultracontractivity for non-symmetric semigroups (see \cite{KS4}).

\medskip

\begin{description}
\item{(A5)}
The transition density function $p^D(t, x, y)$ for $X^D_t$ is
strictly positive in $ D \times D$.
\end{description}

\begin{remark}\label{rk2} {\rm
Even if the L\'evy  process has a smooth and strictly positive
transition density, it is non-trivial to show (A5) (see \cite{B, CZ}
for the case of killed Brownian motions in a domain, \cite{CS1} for
the case of killed symmetric stable processes in a domain and
\cite{V}  for the case of killed non-symmetric stable processes in a
domain). If the L\'evy measure $\nu$ satisfies (A1)(b), the distance
between connected components of $D$ shouldn't be too far away,
otherwise $p^D(t,x,y)$ will be zero there. In section 4, we will
show that for a large class of non-symmetric L\'evy processes, (A5)
is true.}
\end{remark}

In the remainder of this section we always assume that (A1)-(A5) are
in force.

We use $A_D$ and $\wh{A}_D$ to denote the $L^2$ generators of
$\{P^D_t\}$ and   $\{\wh{P}^D_t\}$ respectively. Since for each
$t>0$,  $p^D(t, x, y)$ is bounded in $D\times D$  by (A3),
$\{P^D_t\}$ and   $\{\wh{P}^D_t\}$ are compact operators in $L^2(D,
dx)$. Moreover $p^D(t, x, y)$ is strictly positive in  $ D \times D$
by (A5). Thus it follows from Jentzsch's Theorem (Theorem V.6.6 on
page 337 of \cite{Sc}) and the strong continuity of  $\{P^D_t\}$ and
$\{\wh{P}^D_t\}$ that the common value $\lambda_0:= \sup {\rm Re}
(\sigma(A_D))= \sup {\rm Re} (\sigma(\wh{A}_D)) <0$ is an eigenvalue
of multiplicity 1 for both $A_D$ and $\wh{A}_D$, and that an
eigenfunction $\phi_0$ of $A$ associated with $\lambda_0$ can be
chosen to be strictly positive  a.e. with $\|\phi_0\|_{L^2(D)}=1$
and an eigenfunction $\psi_0$ of $\wh A_D$ associated with
$\lambda_0$ can be chosen to be strictly positive  a.e. with
$\|\psi_0\|_{L^2(D)}=1$. Thus for a.e. $(x,y) \in D \times D$,
\begin{eqnarray}
e^{\lambda_0t} \phi_0 (x) \,=\, \int_D p^D(t,x,z)  \phi_0 (z) dz,
&\quad& -\frac1\lambda_0\, \phi_0 (x) \,=\, \int_D G_D(x,z)  \phi_0
(z) dz,
\label{e:e1}\\
e^{\lambda_0t} \psi_0 (y)\, =\, \int_D \wh p^D(t,y,z)  \psi_0 (z)
dz, &\quad& -\frac1\lambda_0 \,\psi_0 (y) \,=\, \int_D \wh G_D(y,z)
\psi_0 (z) dz. \label{e:e2}
\end{eqnarray}

\medskip

\begin{prop}\label{p:cont}
$\phi_0(x)$ and $\psi_0(x)$ are strictly positive and continuous in
$D$.  Thus  (\ref{e:e1}) and (\ref{e:e2}) are true for every $(x,y)
\in D \times D$.
\end{prop}

\pf By (\ref{e:e1}),
$$
\phi_0 (x) \,=\,e^{-\lambda_0} \int_D p^D(1,x,z)  \phi_0 (z) dz.
$$
Since $p^D(1,x,z)$ is bounded continuous and $D$ is a bounded open
set, the right hand side of the above equation is continuous by
using the dominated convergence theorem and the fact
$\|\phi_0\|_{L^2(D)}=1$. Similarly, $e^{-\lambda_0} \int_D \wh
p^D(1,y,z) \psi_0 (z) dz$ is continuous. Thus there exist continuous
versions of $\phi_0$ and $\psi_0$, and (\ref{e:e1})-(\ref{e:e2}) are
true for every $(x,y) \in D \times D$. Now the strict positivity of
$\phi_0$ and $\psi_0$ follow from the strict positivity of
$p^D(1,\cdot,\cdot)$ and (\ref{e:e1})-(\ref{e:e2}). \qed

\medskip

\begin{defn}
The semigroups $\{P^D_t\}$ and $\{\wh{P}^D_t\}$ are said to be
intrinsic ultracontractive if, for any $t>0$, there exists a
constant $c_t>0$ such that
$$
p^D(t, x, y)\le c_t\phi_0(x)\psi_0(y), \quad \forall (x, y)\in
D\times D.
$$
\end{defn}

\medskip

For results on  intrinsic ultracontractivity for general
non-symmetric semigroups, we refer our readers to Section 2 of
\cite{KS4}.

We will show that the semigroup of any killed non-symmetric L\'evy
process $X^D$ satisfying (A1)-(A5) is intrinsic ultracontractive.
\medskip

\begin{lemma}\label{l:EP}
There exists a constant $c>0$ such that
\begin{equation}\label{e:EP}
\E_x [ \tau_D] \,\le\, c\,\phi_0(x) \quad \mbox{and} \quad \E_y [
\wh \tau_D] \,\le\,c\,\psi_0(y)\qquad
 \forall (x, y)\in D\times D.
\end{equation}
\end{lemma}

\pf  By Lemma \ref{l:main}, there exists a constant $c_1>0$ such
that
$$
\E_x [ \tau_D]\,=\, \int_{B_2} G_D(x,z)dz + \int_{D \setminus B_2}
G_D(x,z)dz  \,\le\, \int_{B_2} G_D(x,z)dz + c_{1}\int_{ B_2}
p^D(t_0, x,z)dz. $$
 Thus by Proposition \ref{p:cont}, we have
\begin{eqnarray*}
&&\int_{B_2} G_D(x,z)dz + c_{1}\int_{ B_2} p^D(t_0, x,z)dz\\
&&\le\, c_2\, \left(
\int_{B_2} G_D(x,z)\phi_0(z)  dz + c_{1}\int_{ B_2} p^D(t_0, x,z)\phi_0(z)dz
\right)\\
&&\le\, c_2\, \left(\int_{D} G_D(x,z)\phi_0(z)  dz +
c_{1}\int_{ D} p^D(t_0, x,z)\phi_0(z)dz
\right)\,=\, c_2\,  \left(-\frac1\lambda_0   + c_{1}  e^{\lambda_0t_0}
\right)
\phi_0 (x)
\end{eqnarray*}
for some positive constant $c_2$.
In the last equality above, we have used (\ref{e:e1}).

Using Lemma \ref{l:main}, Proposition \ref{p:cont} and (\ref{e:e2}),
the second inequality in (\ref{e:EP}) can be proved similarly.
\qed

\medskip

\begin{thm}\label{t:IU}  The semigroups $\{P^D_t\}$ and
$\{\wh P^D_t\}$ are intrinsic
ultracontractive. Moreover, for any $t>0$, there exists a constant
$c_t>0$ such that
\begin{equation}\label{e:IU}
c_t^{-1} \, \phi_0(x)\psi_0(y)
 \, \le \,  p^D(t, x, y)\le \,  c_t \, \phi_0(x)\psi_0(y),
\quad \forall (x, y)\in D\times D.
\end{equation}
\end{thm}

\pf By (A3) and the semigroup property, there exists $c_1(t)>0$ such
that
\begin{eqnarray*}
p^D(t, x, y) &=& \int_D p^D(\frac{t}3, x, z)
\int_D p^D(\frac{t}3, z, w)
p^D(\frac{t}3, w, y) dw dz\\
&\le& c_1(t)  \int_D p^D(\frac{t}3, x, z)
dz\int_D
p^D(\frac{t}3, w, y) dw\\
&=& c_1(t)\, \P_x ( \tau_D > t/3) \, \P_y ( \wh \tau_D > t/3).
\end{eqnarray*}
By applying Chebyshev's inequality
we get
$$
p^D(t, x, y) \,\le\, \frac{ c_1(t) }{9t^2} \,
\E_x [ \tau_D]\,\E_y [ \wh \tau_D].
$$
Thus the intrinsic ultracontractivity is proved by Lemma \ref{l:EP}.

The fact that intrinsic ultracontractivity
implies the lower bound is proved in \cite{KS4}
(Proposition 2.4 in \cite{KS4}). \qed

The following lower bound of $G_D(x,y)$
is an easy corollary of Lemma \ref{l:EP} and
Theorem \ref{t:IU}.

\medskip

\begin{corollary}\label{c:lbG}
There exist constants $c_i>0, i=1,2$ such that
\begin{equation}\label{lbG}
c_1\, \E_x [ \tau_D] \, \E_y [ \wh \tau_D]\, \le\,c_2 \phi_0(x)\psi_0(y)
\le\, G_D(x,y), \qquad (x,y) \in D \times D.
\end{equation}
Moreover, there exists  constant $c_3>0$ such that
$$
c_3^{-1}\,\E_x [ \tau_D]
\,\le\,\phi_0(x)\,\le\,c_3\,\E_x [ \tau_D]
\quad
\mbox{and}
\quad
c_3^{-1}\,\E_x [ \wh \tau_D]
\,\le\,\psi_0(x)\,\le\,c_3\,\E_x [ \wh \tau_D]\qquad
 \forall x\in D.
$$
\end{corollary}

\medskip

Applying Theorem 2.4 in  \cite{KS4}, we have the following.

\begin{thm}\label{c2e}
There exist positive constants $c$ and $\nu$ such that
\begin{equation}\label{h}
\left|\left(e^{-\lambda_0t} \int_D \phi_0(z)\psi_0(z)dz\right)
\frac{p^D(t, x, y)}{\phi_0(x)\psi_0(y)}-1\right|
\le ce^{-\nu t}, \qquad (t, x, y)\in (1,
\infty)\times D\times D.
\end{equation}
\end{thm}

\medskip

We recall the following simple lemma from \cite{KS5}.

\medskip

\begin{lemma}\label{l:mar} (Lemma 5.5 in \cite{KS5})
\begin{description}
\item{(1)}$$\frac{p^D(t, x, y)}{p^D(t, x, z)}\,\ge\,c_1\,
\frac{p^D(t, v, y)}{p^D(t, v, z)},             \quad \forall v,x,y,z
\in D
$$
implies that for every $s >  t  $,
$$
\frac{p^D(s, y, x)}{p^D(s, z, x)}\,\ge\,c_1\, \frac{p^D(t, y,
v)}{p^D(t, z, v)} \quad \mbox{and} \quad\frac{p^D(s, x, y)}{p^D(s,
x, z)}\,\le\,c_1^{-1}\, \frac{p^D(t, v, y)}{p^D(t, v, z)}, \quad
\forall v,x,y,z \in D.
$$

\item{(2)}
$$
\frac{p^D(t, y, x)}{p^D(t, z, x)}\,\ge\,c_2\,
\frac{p^D(t, y, v)}{p^D(t, z, v)}, \quad \forall v,x,y,z \in D
$$implies that for every $s >  t  $,
$$\frac{p^D(s, x, y)}{p^D(s, x, z)}\,\ge\,c_2\,
\frac{p^D(t, v, y)}{p^D(t, v, z)}
\quad
\mbox{and}
\quad
\frac{p^D(s, y, x)}{p^D(s, z, x)}\,\le\,c_2^{-1}\,
\frac{p^D(t, y, v)}{p^D(t, z, v)},\quad  \forall v,x,y,z \in D.
$$
\end{description}
\end{lemma}

\medskip

The parabolic boundary Harnack principle is an easy corollary of
Theorem \ref{t:IU}.

\medskip

\begin{corollary}\label{c:PBHP}
For each positive $u$
there exists $c=c(D,u)>0$ such that
\begin{equation}\label{e:PBHP}
\frac{p^D(t, x, y)}{p^D(t, x, z)}\,\ge\,c\, \frac{p^D(s, v,
y)}{p^D(s, v, z)}, \quad \frac{ p^D(t, y, x)}{ p^D(t, z,
x)}\,\ge\,c\, \frac{ p^D(s, y, v)}{ p^D(s, z, v)}
\end{equation}
for every $s,t \ge u$ and $v,x,y,z \in D$.
\end{corollary}
\pf
By Theorem \ref{t:IU},
both inequalities in (\ref{e:PBHP}) are true for
$s=t=u$.
 Now we apply Lemma \ref{l:mar} (1)-(2) and we get for $ s>u$
\begin{equation}\label{e:L_1}
\frac{p^D(s, y, x)}{p^D(s, z, x)}\,\ge\,c\,
\frac{p^D(u, y, v)}{p^D(u, z, v)},
\quad
\frac{p^D(s, x, y)}{p^D(s, x, z)}\,\le\,c^{-1}\,
\frac{p^D(u, v, y)}{p^D(u, v, z)}, \qquad \forall v,x,y,z \in D
\end{equation}
and
\begin{equation}\label{e:L_2}
\frac{p^D(s, x, y)}{p^D(s, x, z)}\,\ge\,c\,
\frac{p^D(u, v, y)}{p^D(u, v, z)},
\quad
\frac{p^D(s, y, x)}{p^D(s, z, x)}\,\le\,c^{-1}\,
\frac{p^D(u, y, v)}{p^D(u, z, v)},\qquad \forall v,x,y,z \in D.
\end{equation}
Thus both inequalities in (\ref{e:PBHP}) are true for
$s>t=u$. Moreover,
Combining (\ref{e:L_1})-(\ref{e:L_2}),  both inequalities in
(\ref{e:PBHP}) are true for
$t=s>u$ too.  Now applying Lemma \ref{l:mar} (1)-(2) again,
we get our conclusion.
\qed

\medskip

A Borel function $h$ defined on $D$  is said to be
superharmonic
with respect to $X^D$ if
$$
h(x)\,\geq\,  \E_x\left[h(X^D_{\tau_{B}})\right], \qquad x\in B,
$$
for every bounded open set $B$ with $\overline{B}\subset D$.
We use
$SH^+$ to denote families of nonnegative superharmonic functions
of $X^D$.
 For any $h\in SH^+$, we use $\P_x^h$ to denote the
law of the $h$-conditioned process $X^D$
and use $\E_x^h$ to denote the expectation with respect to
$\P_x^h$.
i.e.,
$$
\E_x^h \left[g(X^D_t) \right]\,=\,
\E_x\left[\frac{h(X^D_t)}{h(x)}g(X^D_t) \right]
.$$
Let  $\zeta^h$ be the lifetime of the
$h$-conditioned process $X^D$.

The bound for the lifetime of the
conditioned $X^D$ can be proved using Theorem \ref{c2e}.
It is proved in \cite{KS4} for second order elliptic operators
with smooth coefficients.
Since the proof is similar, we omit the proof here.

\medskip

\begin{thm} (Theorem 3.8 in \cite{KS4})
\begin{description}
\item{(1)}
$$
\sup_{x\in D, h\in SH^+}\E_x^h[\zeta^h]<\infty.
$$
\item{(2)} For any $h\in  SH^+$, we have
$$
\lim_{t\uparrow\infty}e^{-\lambda_0 t}\P_x^h(\zeta^h>t)=
\frac{\phi_0(x)}{h(x)}\int_D\psi_0(y)h(y)dy \Big/\int_D
\phi_0(y)\psi_0(y)dy.
$$
In particular,
$$
\lim_{t\uparrow\infty}\frac1t\log \P_x^h(\zeta^h>t)
=\lambda_0.
$$
\end{description}\end{thm}

\section{Examples}

In this section we collect some examples of L\'evy processes $X$ and
open sets $D$ so that $X^D$ satisfies the assumptions (A1)-(A5).

\medskip

\begin{example}\label{eg1}{\rm
We first recall the definition of non-symmetric strictly
$\alpha$-stable processes. Let  $\alpha\in (0, 2)$ and $d \ge2$. The
process $X$ is said to be strictly $\alpha$-stable if $(X_{at},
\P_0)_{ t \ge 0}$ is equal to $(a^{1/\alpha}X_{t}, \P_0)_{t \ge 0}$
in distribution. Since  $\alpha\in (0, 2)$, $A=0$ and there is a
finite measure $\eta $ on the unit sphere $S=\{  x\in \R^d: |x|=1
\}.$ such that
$$
\nu(U)= \int_S \int_0^\infty 1_U(rz)r^{-(1+\alpha)} dr \eta  (dz)
$$
for every Borel set $U$ in $\R^d$. The measure $\eta$ is called the
spherical part of the L\'{e}vy measure $\nu$. A strictly
$\alpha$-stable process $X$ can be described using its
characteristic function as follows:
\begin{description}
\item{(i)} for $\alpha\in (0, 1)$, a L\'evy process $X$ in $\R^d$ is strictly
$\alpha$-stable if and only if
$$
\E_0\left[e^{iz\cdot
X_1}\right]=\exp\left(\int_S\eta(d\xi)\int^{\infty}_0 (e^{irz\cdot
\xi}-1)r^{-(1+\alpha)}dr\right);
$$
\item{(ii)} for $\alpha=1$, a L\'evy process $X$ in $\R^d$ is strictly
$\alpha$-stable if and only if
$$
\E_0\left[e^{iz\cdot
X_1}\right]=\exp\left(\int_S\eta(d\xi)\int^{\infty}_0 (e^{irz\cdot
\xi}-1-irz\cdot\xi1_{(0, 1]}(r))r^{-2}dr+iz\cdot \gamma\right)
$$
for some $\gamma\in \R^d$ and $\int_S\xi\eta(d\xi)=0$;
\item{(iii)} for $\alpha\in (1, 2)$, a L\'evy process $X$ in $\R^d$
is strictly $\alpha$-stable if and only if
$$
\E_0\left[e^{iz\cdot
X_1}\right]=\exp\left(\int_S\eta(d\xi)\int^{\infty}_0 (e^{irz\cdot
\xi}-1-irz\cdot\xi)r^{-(1+\alpha)}dr\right).
$$
\end{description}

Suppose that $X=(X_t, \P_x)$ is a strictly $\alpha$-stable process
with the spherical part $\eta$ of its L\'{e}vy measure satisfying
the following assumption:
 there exist $\varphi: S \to (0, \infty)$ and $\kappa > 0$ such that
\begin{equation}\label{e:slevy}
\varphi= \frac{d \eta }{d \sigma} \quad \mbox{ and } \quad \kappa \,
\le \, \varphi (z) \, \le \,  \kappa^{-1}, \quad \forall z \in S,
\end{equation}
where $\sigma$ is the surface measure on $S$.
Thus the L\'{e}vy measure $\nu$ has a density
$f(x)=\varphi (x/|x|) |x|^{-(d+ \alpha)}$ with respect to the
$d$-dimensional Lebesgue measure, and
\begin{equation}\label{e:Levy_e}
\kappa \,|x|^{-(d+\alpha)}\, \le\, f(x)\, \le\,
\kappa^{-1}\, |x|^{-(d+\alpha)},
\qquad x \in \R^d.
\end{equation}
Thus it is easy to see that (A1)(a) is true with
$L(x)=|x|^{d+\alpha}\varphi (x/|x|)^{-1}$.

The process $X$ has a
jointly continuous and strictly positive transition
density function $p(t,x,y)=p(t,x-y)$
and there exists $c > 0$ such that
\begin{equation}\label{e:p_up}
p(t,x,y) \,\le\,  c\, t^{-\frac{d}{\alpha}},\qquad \forall
(t,x,y) \in (0,\infty)
\times \R^d \times \R^d
\end{equation}
(see (2.6) in \cite{V}).
Moreover, for any $\gamma>0$, there exists $c>0$ such that
\begin{equation}\label{e:p_upod}
p(t,x,y) \,\le\,  c\, t, \qquad |x-y|\ge \gamma, \,t>0.
\end{equation}
(see (2.5) in \cite{V}). Using the facts above, one can follow
routine arguments (see, for instance, the proof of Theorem 2.4 in
\cite{CZ}) to show that, for every open subset $D$,  the killed
process $X^D$ has a transition density $p^D(t, x, y)$ such that, for
any $t>0$, $p^D(t, x, y)$ is jointly continuous on $D\times D$. It
follows from Theorem 3.2 in \cite{V} that $p^D(t, x, y)$ is strictly
positive on $(0, \infty)\times D\times D$ when $D$ is connected. Now
we are going to show that $p^D(t, x, y)$ is strictly positive on
$(0, \infty)\times D\times D$ when $D$ is not connected. It is
enough to show that for any two connected components $D_1$ and $D_2$
of $D$, $p^D(t, x, y)$ strictly positive on $(0, \infty)\times
D_1\times D_2$. It follows from Lemma \ref{l:main} that for any
$x\in D_1$ and any ball $B(x_0, r)$ with $\overline{B(x_0,
r)}\subset D_2$, there exist constants $t_0>0$ and $c>0$ such that
$$
\P_x(X_t\in B(x_0, r), t<\tau_D)\,\ge\, c\int_{D\setminus B(x_0,
r)}G_D(x, y)dy \,\ge\, c\int_{D_1}G_{D_1}(x, y)dy>0
$$
whenever $t\le t_0$. This implies that, for $t\le t_0$ and $x\in
D_1$, $p^D(t, x, \cdot)$ is strictly positive almost everywhere on
$D_2$. By working with the dual process we get that for $t\le t_0$
and $y\in D_2$, $p^D(t, \cdot, y)$ is strictly positive almost
everywhere on $D_1$. Combining these with the semigroup property we
get that $p^D(t, x, y)$ is strictly positive everywhere on $(0,
\infty)\times D_1\times D_2$. Thus in this case (A2), (A3) and (A5)
 are valid for any bounded open subset $D$ as well. }
\end{example}

\medskip

\begin{example}\label{eg2}{\rm
Assume that $X$ is a non-symmetric strictly $\alpha$-stable
processes from the previous example and we will use the notations
from the previous example. A L\'{e}vy process $Y$ in $\R^d$ is
called truncated (non-symmetric) strictly $\alpha$-stable process if
\begin{description}
\item{(i)} when $\alpha\in (0, 1)$,
$$
\E_0\left[e^{iz\cdot Y_1}\right]=\exp\left(\int_S\eta(d\xi)\int^1_0
(e^{irz\cdot \xi}-1)r^{-(1+\alpha)}dr\right);
$$
\item{(ii)} when $\alpha=1$,
$$
\E_0\left[e^{iz\cdot Y_1}\right]=\exp\left(\int_S\eta(d\xi)\int^1_0
(e^{irz\cdot \xi}-1-irz\cdot\xi)r^{-2}dr+iz\cdot \gamma\right)
$$
for some $\gamma\in \R^d$ and $\int_S\xi\eta(d\xi)=0$;
\item{(iii)} when $\alpha\in (1, 2)$,
$$
\E_0\left[e^{iz\cdot Y_1}\right]=\exp\left(\int_S\eta(d\xi)\int^1_0
(e^{irz\cdot \xi}-1-irz\cdot\xi)r^{-(1+\alpha)}dr\right).
$$
\end{description}
We also assume that $\eta$ satisfies \eqref{e:slevy}. Then the
L\'{e}vy density $g(x)$ for $Y$ is
\begin{equation}\label{e:Levy_t}
g(x)\,:=\,\varphi (x/|x|) |x|^{-(d+ \alpha)} 1_{\{|x| <1\}}
\end{equation}
and (A1)(b) is satisfied. In the case when $Y$ is rotationally
invariant, it has been studied recently by the authors \cite{KS}.
\eqref{e:slevy} implies that the characteristic function of $Y_t$ is
integrable. Thus the process $Y$ has a bounded and continuous
density $q(t, x, y)$ (cf. \cite{Sa}). Let
\begin{equation}\label{e:pl}
h (x) \,:= \,f(x)-g(x)\,=\,
\varphi (x/|x|) |x|^{-(d+ \alpha)} 1_{\{|x| \ge 1\}}.
\end{equation}
 Note that $\lambda:=\int_{\R^d} h (x)  dx <\infty$.
Thus we can write $X_t =Y_t +Z_t$ where
$Z_t$
 is a compound Poisson process with the L\'evy density $h(x)$,
independent of $Y_t$.
Let
$$
T\,:=\,\inf\{t\ge 0: V_t \not=0\}.
$$
$T$ is an exponential random variable with intensity $\lambda$.
Moreover,
 $Y_t=X_t$ for $t < T$  and $\{t < \tau^Y_D,\, t<T\}
= \{t < \tau^X_D,\, t<T\}$ where  $\tau^X_D:=\inf\{t>0: \, X_t\notin D\}$ and
 $\tau^Y_D:=\inf\{t>0: \, Y_t\notin D\}$.
Thus, since $Y$ and $T$ are independent, for every open subsets $U$
and $D$ with $U \subset D$ we have
\begin{eqnarray}
\P(Y_t^D \in U\,|\,Y_0=x)\P(T>t)&=&
\P(Y_t \in U,\, t < \tau^Y_D,\, t<T\,|\,Y_0=x)\nonumber\\
&=&
\P(X_t \in U,\,t < \tau^X_D,\, t<T\, |\,X_0=x)\label{e:pp1}\\
&\le&
\P(X^D_t \in U,\,t < \tau^X_D\, |\,X_0=x)\label{e:pp2}
\end{eqnarray}
One can find a similar argument for symmetric L\'evy processes in
the proof of Lemma 2.5 in \cite{GR}.

From (\ref{e:pp2}) with $D=\R^d$,
 we have
\begin{equation}\label{e:p_st}
q(t,x,y) \, \le \, e^{\lambda t} p(t,x,y), \qquad \forall (t,x,y)
\in (0, \infty)\times \R^d \times\R^d.
\end{equation}
Combining  (\ref{e:p_up}), (\ref{e:p_upod}) and (\ref{e:p_st}), we get
\begin{equation}\label{e:p_upt}
q(t,x,y) \,\le\,  c\,e^{\lambda t} t^{-\frac{d}{\alpha}},\qquad \forall
(t,x,y) \in (0,\infty)
\times \R^d \times \R^d
\end{equation}
and
\begin{equation}\label{e:p_upodt}
q(t,x,y) \,\le\,  c\, t\, e^{\lambda t}, \qquad |x-y|\ge \gamma,
\,t>0.
\end{equation}
Using the facts above, one can follow routine arguments (see,
for instance, the proof of Theorem 2.4 in \cite{CZ}) to show
that,
for every open subset $D$,  the killed process $Y^D$ has a
transition density $q^D(t, x, y)$ such that, for any $t>0$,
$q^D(t, x, y)$ is jointly continuous on $D\times D$.

Due to our assumption on the L\'evy measure of $Y$, $q^D(t, x, y)$
may not be strictly positive without further assumption on the open
set $D$. Now we are going to show that, when $D$ is a bounded
roughly connected open set, $q^D(t, x, y)$ is strictly positive on
$(0, \infty)\times D\times D$. Thus in this case, (A2), (A3) and
(A5) are satisfied.}
\end{example}

\begin{defn}\label{d:rc}
We say that an open set $D$ in $\R^d$ is roughly connected if
 for every $x, y \in D$, there exist distinct connected components
$U_1 \cdots, U_m$ of $D$ such that
 $ x \in U_1$, $y \in U_m$ and dist$(U_k, U_{k+1}) <1$ for $1 \le k \le m-1$.
\end{defn}

\begin{prop}\label{p:spdp}
For every bounded roughly connected open set $D$, the transition
density function $q^D(t,x,y)$ for $Y$ in $D$ is strictly positive in
$(t,x,y) \in (0, \infty)\times D \times D$.
\end{prop}

\pf We prove the proposition in several steps.
\begin{description}
\item{(1)} We first assume that $\mbox{diam}(D) < 1$.
Fix $t>0$.
We recall from (\ref{e:pp1}) that for every non-empty open set $U \subset D$
$$
\P_x(Y_{t/2}^D \in U)\,=\,e^{\lambda t/2}
\P(X_{t/2} \in U,\,t/2 < \tau^X_D,\, t/2<T\, |\,X_0=x).
$$
Note that by (\ref{e:pl}), we know $Z_t$ makes jumps with sizes great than
or equal to $1$ only. Thus, since $\mbox{diam}(D) < 1$,
$\{t/2 < \tau^X_D,\, t/2<T\}= \{t/2 < \tau^X_D\}$, which implies that
$$
\int_{U} q^D(t/2,x,y) dy\,=\,
\P_x(Y_{t/2}^D \in U)\,=\,e^{\lambda t/2}\P_x(X_{t/2}^D \in U) >0.
$$
Thus for each $x \in D$, $q^D(t/2,x,y)>0$ for a.e. $y \in D$.
Similarly,
$$
\int_{U} q^D(t/2,x,y) dx\,=\, \P_y(\wh Y_{t/2}^D \in
U)\,=\,e^{\lambda t/2}\P_y(\wh X_{t/2}^D \in U) >0.
$$
Thus, for each $y \in D$, $q^D(t/2,x,y)>0$ for a.e. $x \in D$.
Therefore the semigroup property implies that
$$
 q^D(t,x,y)\,=\, \int_D q^D(t/2,x,z) \, q^D(t/2,z,y)\, dz $$
is strictly positive for $(x,y) \in D \times D$ in this case.
\item{(2)} Now we assume that $D$ is connected.
If $x, y \in D \cap B(x_0, r)$ where $x_0 \in D$ and $r<1/2$, then by (1)
\begin{equation}\label{e:spdp1}
 q^D (t,x,y) \,\ge\, q^{D \cap B(x_0, 1/2)} (t,x,y) \, > \, 0.
\end{equation}
Thus by the semigroup property and (\ref{e:spdp1}), for $y \in B(x, 1/2)$,
\begin{eqnarray}
 q^D (t,x,y) &=& \int_D q^D (t/2,x,z) \,  q^D (t/2,z,y) \, dz \nonumber \\
&\ge & \int_{D\cap B(x, 1/2)}  q^D(t/2,x,z) \,  q^D (t/2,z,y) \, dz
\, >\, 0.\label{e:spdp2}
\end{eqnarray}
Using this and a simple chain argument one can easily show that
$q^D(t, x, y)$ is strictly positive on $(0, \infty)\times D\times D$
in this case.

\item{(3)}  Finally we deal with the general case that $D$ is a roughly connected open set;
Fix $x, y \in D$.  There exist distinct connected components $U_1,
\cdots, U_m$ of  $D$ and $\eps>0$ such that
 $ x \in U_1$, $y \in U_m$
and dist$(U_k, U_{k+1}) <1 -4 \eps  $ for $1 \le k \le m-1$. Choose
points $x^1_{k}, x^2_{k} \in U_k$ and $\delta^1_{k}, \delta^2_{k}
<\eps$ where
 $1 \le k \le m$ such that $x=x^1_{1}$, $y=x^2_{m}$,
$|x^2_{k}-x^1_{k+1}|<1 -2 \eps$ and
$$
V_{k,k+1}\,:= \,
 B(x^2_{k},\delta^2_{k})\cup B(x^1_{k+1}, \delta^1_{k+1}) \, \subset\, U_k \cup U_{k+1},
$$
for  $1 \le k \le m-1$. Let $t_m:=t/(2m-1)$.
Now by the semigroup property
\begin{eqnarray*}
&& q^D (t,x,y)\\
&&=
\int_D \cdots \int_D  q^D (t_m,x^1_{1},y^2_{1} )
q^D (t_m,y^2_{1},y^1_{2} )\cdots
  q^D (t_m,y^1_{k},y^2_{k} )\\
  &&~~~ \times
q^D (t_m,y^2_{k},y^1_{k+1} )\cdots
q^D (t_m,y^2_{m-1},y^1_{m} )  q^D (t_m,y^1_{m},x^2_{m} )
dy^2_1 dy^1_2\cdots  dy^2_{m-1}dy^1_m\\
&&\ge \int_{U_1}  q^D (t_m,x^1_{1},y^2_{1} )
 \int_{ V_{1,2}   }
q^D (t_m,y^2_{1},y^1_{2} )\cdots \int_{U_k}   q^D (t_m,y^1_{k},y^2_{k} )  \\
  &&~~~ \times
  \int_{V_{k,k+1}}   q^D (t_m,y^2_{k},y^1_{k+1} )\cdots
\int_{V_{m-1,m}}  q^D (t_m,y^2_{m-1},y^1_{m} ) \\
  &&~~~ \times
\int_{U_m}  q^D (t_m,y^1_{m},x^2_{m} )   dy^2_1 dy^1_2\cdots
dy^2_{m-1} dy^1_m,
\end{eqnarray*}
which is strictly positive in $D \times D$ by (1)-(2).
 \end{description}
\qed

\medskip

\begin{example}\label{eg3}{\rm Suppose that $X$ is a strictly $\alpha$-stable
process in $\R^d$ satisfying all the assumptions in Example
\ref{eg1}, that $B$ is a Brownian motion in $\R^d$ and that $X$ and
$B$ are independent. Then the process $Z$ defined by $Z_t=B_t+X_t$
is also a L\'evy process and it obviously satisfies (A1)(a). The
transition density $q(t, x, y)$ of $Z$ is given by the convolution
of the transition densities of $B$ and $X$. Using this, the explicit
formula for the transition density of $B$, and \eqref{e:p_up} and
\eqref{e:p_upod} for the transition density of $X$, we can easily
show that there exists $c > 0$ such that
\begin{equation}\label{e:p_upm1}
q(t,x,y) \,\le\,  c\, t^{-\frac{d}{\alpha}},\qquad \forall (t,x,y)
\in (0,\infty) \times \R^d \times \R^d.
\end{equation}
Moreover, for any $\gamma>0$, there exists $c>0$ such that
\begin{equation}\label{e:p_upodm1}
q(t,x,y) \,\le\,  c t, \qquad |x-y|\ge \gamma, \,t>0.
\end{equation}
Using the facts above, one can follow routine arguments (see, for
instance, the proof of Theorem 2.4 in \cite{CZ}) to show that, for
every open subset $D$,  the killed process $Z^D$ has a transition
density $q^D(t, x, y)$ such that, for any $t>0$, $q^D(t, x, y)$ is
jointly continuous on $D\times D$. It follows from the lemma below
that $q^D(t, x, y)$ is strictly positive on $(0, \infty)\times
D\times D$. Thus for any bounded open subset $D$, $Z^D$ satisfied
(A2), (A3) and (A5).}
\end{example}

\begin{lemma}\label{newlemma}
Suppose that $Z$ is the process in Example \ref{eg3} and that $q(t,
x, y)$ is the transition density of $Z$. Then for any bounded open
set $D$ in $\R^d$, the transition density $q^D(t, x, y)$ of $Z^D$ is
strictly positive on $(0, \infty)\times D\times D$.
\end{lemma}

\pf For any bounded domain $V$ and bounded open set $U$, let
$p^V_1(t, x, y)$ and $p^U_2(t, x, y)$ be the density of the killed
Brownian motion $B^V$ and the killed strictly $\alpha$-stable
process $X^U$ respectively. Note that the above densities are
strictly positive.

Without loss of generality we assume $B_0=0$. For $x \in D$, let
$\delta_x$ be a positive constant with $B(x, 2\delta_x) \subset D$.
We will show that for every $B(x_0, \eps) \subset D$, $\P_x(Z^D_t
\in B(x_0, \eps))>0$.

Choose $\delta=\delta(x_0, x, \eps) < \delta_x$ such that $B(x_0,
\eps)\subset B(x_0, \eps + \delta) \subset D $ and let $U:=B(x_0,
\eps) \cup B(x, \delta_x)$. Then
\begin{eqnarray*}
\P_x\left(Z^D_t \in B(x_0, \eps)\right)&=&\P_x\left(B_t+X_t \in
B(x_0, \eps), \,
\tau^Z_D >t\right)\\
&\ge &\P_x\left(B_t+X_t \in B(x_0, \eps), \, \tau^B_{B(0, \delta)}
>t, \,
 \tau^X_{U} >t \right)\\
&=&\int_{B(0, \delta)} \int_{B(x_0, \eps)} p^{U+y}_2(t, x+y, z)
p^{B(0, \delta)}_1(t, 0, y) dzdy,
\end{eqnarray*}
which is strictly positive. This implies that, for $t < \infty$ and
$x\in D$, $q^D(t, x, \cdot)$ is strictly positive almost everywhere
on $D$. By working with the dual process we get that for $t< \infty$
and $y\in D$, $q^D(t, \cdot, y)$ is strictly positive almost
everywhere on $D$. Combining these with the semigroup property we
get that $q^D(t, x, y)$ is strictly positive everywhere on $(0,
\infty)\times D\times D$.
 \qed

\medskip

\begin{example}\label{eg4}{\rm Suppose that $Y$ is a truncated
strictly $\alpha$-stable process in $\R^d$ satisfying all the
assumptions in Example \ref{eg2}, that $B$ is a Brownian motion in
$\R^d$ and that $Y$ and $B$ are independent. Then the process $Z$
defined by $Z_t=B_t+Y_t$ is also a L\'evy process and it obviously
satisfies (A1)(b). The transition density $k(t, x, y)$ of $Z$ is
given by the convolution of the transition densities of $B$ and $Y$.
Using this, the explicit formula for the transition density of $B$,
and \eqref{e:p_upt} and \eqref{e:p_upodt} for the transition density
of $Y$, we can easily show that there exist $c > 0$ and $\lambda>0$
such that
\begin{equation}\label{e:p_upm2}
k(t,x,y) \,\le\,  c\, e^{\lambda t}\,t^{-\frac{d}{\alpha}},\qquad
\forall (t,x,y) \in (0,\infty) \times \R^d \times \R^d.
\end{equation}
Moreover, for any $\gamma>0$, there exist $c>0$ and $\lambda>0$ such
that
\begin{equation}\label{e:p_upodm2}
k(t,x,y) \,\le\,  c\, t\, e^{\lambda t}, \qquad |x-y|\ge \gamma,
\,t>0.
\end{equation}
Using the facts above, one can follow routine arguments (see, for
instance, the proof of Theorem 2.4 in \cite{CZ}) to show that, for
every open subset $D$,  the killed process $Z^D$ has a transition
density $Z^D(t, x, y)$ such that, for any $t>0$, $k^D(t, x, y)$ is
jointly continuous on $D\times D$. It follows from the lemma below
that for every bounded roughly connected open set $D$, $k^D(t, x,
y)$ is strictly positive on $(0, \infty)\times D\times D$. Thus for
any bounded roughly connected open subset $D$, $Z^D$ satisfied (A2),
(A3) and (A5).}
\end{example}

\begin{lemma}\label{newlemma2}
Suppose that $Z$ is the process in Example \ref{eg4} and that $k(t,
x, y)$ is the transition density of $Z$. Then for any bounded
roughly connected open set $D$ in $\R^d$, the transition density
$k^D(t, x, y)$ of $Z^D$ is strictly positive on $(0, \infty)\times
D\times D$.
\end{lemma}

\pf For any bounded domain $V$ and bounded open set $U$, let
$p^V_1(t, x, y)$ and $p^U_2(t, x, y)$ be the density of the killed
Brownian motion $B^V$ and the killed truncated $\alpha$-stable
process $Y^U$ respectively. Note that if diam$(U)<1$, the above
densities are strictly positive (Proposition \ref{p:spdp}). Thus
through the same argument in the proof of Lemma \ref{newlemma}, we
have that for every open subset $D$ with diam$(D)<1$, $k^D(t,x,y)$
is strictly positive. Now following the step (2)-(3) in the proof of
Proposition \ref{p:spdp}, we conclude that for any bounded roughly
connected open set $D$, $k^D(t, x, y)$ is strictly positive on $(0,
\infty)\times D\times D$. \qed

\medskip

We list here more examples of L\'evy processes $X$ and open sets $D$
so that $X^D$ satisfies the assumptions (A1)-(A5) without giving
proofs. One can prove them easily using arguments similar to those
in the previous examples and induction.

If $X^{(j)}, j=1,\dots, n,$ are independent strictly
$\alpha_j$-stable processes satisfying the assumptions of Example
\ref{eg1}. Then the process $X$ defined by $X_t=X^{(1)}_t+
\cdot+X^{(n)}_t$ is a L\'evy process satisfying (A1)(a). For any
bounded open subset of $\R^d$, the killed process $X^D$ satisfies
(A2), (A3) and (A5). Similarly, if $X^{(j)}, j=1,\dots, n,$ are
independent truncated strictly $\alpha_j$-stable processes all
satisfying the assumptions of Example \ref{eg1}. Then the process
$X$ defined by $X_t=X^{(1)}_t+ \cdot+X^{(n)}_t$ is a L\'evy process
satisfying (A1)(b). For any bounded roughly connected open subset of
$\R^d$, the killed process $X^D$ satisfies (A2), (A3) and (A5). Of
course, one can combine the examples above with Brownian motion to
get more examples.

\medskip
\vspace{.1in}
{\bf Acknowledgment}: This paper was finalized while
the first named author was visiting KIAS in Seoul, the hospitality
of which is acknowledged.

 \vspace{.5in}
\begin{singlespace}

\end{singlespace}

\end{doublespace}


\small
\begin{thebibliography}{99}


\bibitem{Ba} R.~Ba\~nuelos,
Intrinsic ultracontractivity and eigenfunction estimates for
Schr\"odinger operators. {\em J. Funct. Anal.}
{\bf 100} (1991), 181-206.

\bibitem{B}  R. F. Bass,
{\it  Probabilistic Techniques in Analysis}.
 Springer-Verlag, 1995.


\bibitem{BB}
R. F. Bass and K. Burdzy,
Lifetimes of conditioned diffusions.
{\em Probab. Theory Related Fields} {\bf 91} (1992), 405--443.


\bibitem{Be}
J. Bertoin,
{\em L\'evy processes,}
Cambridge University Press, Cambridge, 1996.



\bibitem{BG}  R.~M. Blumenthal and  R.~K. Getoor,
{\it  Markov Processes and Potential Theory}.
Academic Press, 1968.



\bibitem{CS1}
Z.-Q.~Chen and R.~Song, Intrinsic ultracontractivity and conditional gauge
for symmetric stable processes.
{\em J. Funct. Anal.} {\bf 150} (1997), 204--239.

\bibitem{CS2}
Z.-Q.~Chen and R.~Song,
Intrinsic ultracontractivity, conditional lifetimes and
conditional gauge for symmetric stable processes on rough domains.
{\em Illinois J. Math.} {\bf 44(1)} (2000), 138--160.

\bibitem{C2} K. L. Chung,
Greenian bounds for Markov processes,
{\em Potential Anal. } {\bf 1(1)} (1992), 83--92.



\bibitem{CW} K.~L. Chung and J.~B. Walsh,
{\it Markov processes, Brownian motion, and time symmetry}.
Springer, New York, 2005


\bibitem{CZ}
K. L. Chung and Z. X. Zhao,
{\em From Brownian motion to Schr\"odinger's equation},
Springer-Verlag, Berlin, 1995.


\bibitem{DS}
E.~B. Davies and B. Simon.
Ultracontractivity and the heat kernel for Schr\"{o}dinger
operators and Dirichlet Laplacians.
{\em J. Funct. Anal.} {\bf 59} (1984), 335-395.


\bibitem{G} T. Grzywny, Intrinsic ultracontractivity for L\'evy
processes,
Preprint, 2006.

\bibitem{GR}
T. Grzywny and M. Ryznar, Estimates of Green function for some
perturbations of fractional Laplacian,
preprint, 2006


\bibitem{IW}
N. Ikeda and S. Watanabe,
On some relations between the harmonic measure and the
L\'{e}vy measure for a certain class of Markov processes, {\em
J. Math. Kyoto Univ.} {\bf 2} (1962), 79--95.


\bibitem{JK} D. S. Jerison and C. E. Kenig,
Boundary behavior of harmonic functions in non-tangentially
accessible domains. {\it Adv. Math.}, {\bf 46} (1982), 80--147.

\bibitem{KS} P. Kim and R. Song, Potential theory of truncated
stable processes. {\it Math. Z.} (to appear), 2006.

\bibitem{KS4} P. Kim and R. Song, Intrinsic ultracontractivity of
non-symmetric diffusion semigroups in bounded domains,
Preprint, 2006.

\bibitem{KS5} P. Kim and R. Song,
Intrinsic ultracontractivity of non-symmetric diffusions with measure-valued
drifts and potentials, Preprint, 2006.



\bibitem{Ku}
T. Kulczycki,
Intrinsic ultracontractivity for symmetric stable processes,
{\it Bull. Polish Acad. Sci. Math.} {\bf 46(3)} (1998), 325--334.


\bibitem{MV} O. Martio and M. Vuorinen,
Whitney cubes, $p$-capacity, and Minkowski content. {\it Exposition.
Math.}, {\bf 5(1)} (1987), 17--40.



\bibitem{Sa}
K. Sato, {\em
L\'evy processes and infinitely divisible distributions},
Cambridge University Press, Cambridge, 1999.



\bibitem{Sc} H. H. Schaefer,
{\em Banach lattices and positive operators},
Springer-Verlag, New York, 1974.


\bibitem{SW}
R. Song and J. Wu, Boundary Harnack principle for symmetric stable
processes, {\it J. Funct. Anal.} {168(2)} (1999), 403-427.


\bibitem{V}
Zoran, Vondra\v cek,
Basic potential theory of certain nonsymmetric strictly
$\alpha$-stable processes.
{\em Glas. Mat. Ser. III} {\bf 37(57)} (2002), 211--233.





\end{thebibliography}
\end{document}